\documentclass[11pt]{article}
\usepackage{amssymb,amsmath,latexsym,amscd,amsfonts}
\usepackage{graphics}
\usepackage{epsfig}

\newtheorem{thm}{Theorem}[section]

\newtheorem{definition}{Definition}[section]

\newtheorem{lem}{Lemma}[section]

\newtheorem{cor}{Corollary}[section]

\newenvironment{proof}{{\noindent{\bf Proof: }}}{$\hfill\Box$}

\newcommand{\ignore}[1]{}

\newcommand{\onote}[1]{}
\newcommand{\dfinalnote}[1]{}
\newcommand{\znote}[1]{}

\def\hpic #1 #2 {\mbox{$\begin{array}[c]{l}
\epsfig{file=#1,height=#2} \end{array}$}}
\def\vpic #1 #2 {\mbox{$\begin{array}[c]{l}
\epsfig{file=#1,width=#2}\end{array}$}}

\begin{document}
\title{ Enumerating the Classes of Local Equivalency in Graphs}

\author{{Mohsen Bahramgiri} \thanks{Massachusetts Institute of Technology,
Mathematics department \& Computer Science and Artificial
Intelligence Laboratory,\ m\underline{ }bahram@mit.edu.}
 \and Salman Beigi\thanks{Massachusetts Institute of Technology,
Mathematics department,\ salman@math.mit.edu.} }

\date{February 09, 2007}

\maketitle{}

\begin{abstract} There are local operators on (labeled) graphs $G$ with labels
$(g_{ij})$ coming from a finite field. If the filed is binary, in
other words, if the graph is ordinary, the operation is just the
local complementation. That is, to choose a vertex and complement
the subgraph induced by its neighbors. But, in the general case,
there are two different types of operators. The first type is the
following. Let $v$ be a vertex of the graph and $a\in
\mathbf{F}_q$, the finite field of $q$ elements. The operator is
to obtain a graph with labels $g'_{ij}=g_{ij}+ag_{vi}g_{vj}$. For
the second type of operators, let $0\neq b\in \mathbf{F}_q$ and
the resulted graph is a graph with labels $g''_{vi}=bg_{vi}$ and
$g''_{ij}=g_{ij}$, for $i,j$ unequal to $v$.

The local complementation operator (binary case) has appeared in
combinatorial theory, and its properties have studied in the
literature, \cite{bouchalg, bouchtree, bouchet}. Recently, a
profound relation between local operators on graphs and quantum
stabilizer codes has been found \cite{moor1, us}, and it has
become a natural question to recognize equivalency classes under
these operators. In the present article, we show that the number
of graphs locally equivalent to a given graph is at most
$q^{2n+1}$, and consequently,  the number of classes of local
equivalency is $q^{\frac{n^2}{2}-o(n)}$.
\end{abstract}

\section{Introduction}

A {\it labeled graph} is a graph with labeled edges, with labels
coming from a (finite) field. This covers the ordinary (simple)
graphs, when one restricts the field to be the binary field,
$\mathbf{F}_2$. For simplicity, we discuss the binary case
separately to make the notion more clear. In the binary case,
consider the following operator, called {\it local
complementation}. Choose a vertex, and replace the graph induced
by the neighbors of this vertex by its complement. Two graphs are
called {\it locally equivalent} if one can be obtained from the
other by applying some local operations described earlier.

In general, when the field is not binary, there are two
independent types of operators involved. The first one is just the
generalized version of the operation in binary case. Let the
graph $G$ be labeled with labels forming a symmetric matrix
$G=(g_{ij})$ with zero diagonal, over $\mathbf{F}_q$ where $q$ is
a power of a prime number, and $\mathbf{F}_q$ is the field with
$q$ elements. Let $v$ be a vertex of this graph, and $a\in
\mathbf{F}_q$. We define the first type of operators in the
following way.  $G*_{a}v$ is defined to be a graph with labels
$G'=(g'_{ij})$ such that $g'_{ij}=g_{ij}+a \cdot g_{vi}g_{vj}.$
The second type of operators is multiplying the edges of some
vertex by a non-zero number, $b \in \mathbf{F}_q$. In other words,
$G\circ_{b}v$ is the graph with labels $G''=(g''_{ij})$, where
$g''_{vi}=bg_{vi}$ and $g''_{ij}=g_{ij}$ for $i, j$ unequal to
$v$. Similar to the previous situation, two graphs are called
locally equivalent if one can obtain one of them by applying the
operators $*$ and $\circ$ on the other one.

Studying and investigating the local equivalency of graphs has
become a natural problem in quantum computing, and playing a
significant role especially in {\it error-correcting codes}, due
to the recent work of \cite{knill}, \cite{eff}, \cite{moor1},
\cite{us} and \cite{marc} . Namely, in the quantum computing
setting, some states, called {\it stabilizer states}, have a
description as the common eigenvector of a subgroup of the {\it
Pauli group}. Using stabilizer states, we may be able to create
more preferable {\it quantum codes}, due to the property that the
obtained codes, have relatively shorter description to handle the
process. On the other hands, {\it graph states}, an important
subset of stabilizer states, are defined based on graphs with
labels in a finite field. Combining the theory of
error-correcting codes and the tools in generalized graph theory,
leads us to describe and investigate the properties of graph
states more and more deeply.

Some stabilizer states may have similar properties. In fact, we
can obtain one of the stabilizer states from the other by applying
elements of {\it local Clifford group}. If two states are
equivalent under local Clifford group, they present similar
properties in quantum computing. The key point is that, any
stabilizer state is equivalent to a graph state under the local
Clifford group, and consequently, we may just consider the graph
states. On the other hand, some of the graph states are
equivalent under the local Clifford group. More precisely, shown
in \cite{moor1} and \cite{us}, two graph states are equivalent
under the local Clifford group if their associated graphs are
locally equivalent by the local operations described earlier. The
properties of locally equivalent graphs have been deeply studied
in the recent works, and an efficient algorithm to determine
whether two graphs are locally equivalent or not is given. This
algorithm for the binary case can be found in \cite{bouchalg},
and for the general case in \cite{eff}.

The purpose of this article is solving a significantly crucial
problem in studying the local operations: enumerating the graphs
locally equivalent to a given one, as well as enumerating the
equivalency classes.

\subsection{Main results}

The main results proven in the present paper are the followings.

First, The number of graphs locally equivalent to a given one is
at most $(q-1)(q^2-1)^n$, which is bounded above by $q^{2n+1}$,
$n$ being the number of vertices of the graph.

Second, $\mathcal{C}(n)$, the number of classes of local
equivalency of graphs with $n$ vertices satisfies:
$$q^{\frac{n^2}{2}-\frac{5n}{2}-1} \leq \mathcal{C}(n) \leq q^{\frac{n^2}{2}-\frac{n}{2}}. $$
In other words, $$\mathcal{C}(n) = q^{\frac{n^2}{2}-O(n)}.$$ In
particular, for the usual (binary) graphs, i.e. when $q=2$, the
number of graphs locally equivalent to a graph is at most $3^n$
and the number of classes of local equivalency is
$$\mathcal{C}(n) = 2^{\frac{n^2}{2}-O(n)}.$$

\subsection{Structure of the paper}

This paper is organized as follows; In section 2, we introduce
the geometrically known concept of {\it isotropic systems} and
also a relation between isotropic systems and graphs. In fact, we
correspond to every graph an isotropic system, and say that graph
is a {\it graphic presentation} for the isotropic system. Then, we
show that every isotropic system has a graphic presentation.
This, somehow, says that the properties of graphs and isotropic
systems are involved.

In section 3, after introducing the definitions of local
operators, in theorem \ref{loeq}, we translate local equivalency
into an algebraic equation, which is significantly helpful
throughout this article. We will then prove that two fundamental
graphs for an isotropic system are (up to a constant) locally
equivalent.

Eulerian vectors, which, roughly speaking, are orthogonal vectors
to an isotropic system, are introduced rigorously in section 4.
The number of Eulerian vectors for a given isotropic system, say
$\mathcal{L}$, is denoted by $\epsilon(\mathcal{L})$, and it is
shown that, $\epsilon(\mathcal{L})\geq 1$ for every
$\mathcal{L}$. The notion of {\it switching property} is
introduced in this section, and using its power, the exact number
of graphic presentations is given in terms of
$\epsilon(\mathcal{L})$.

In section 5, we introduce the notion of {\it index of an
isotropic system}, denoted by $\lambda(\mathcal{L})$, and then we
estimate it, from above as well as below, by the terms containing
the dimension of the {\it bineighborhood} space, introduced and
studied in this section.

Section 6 is dedicated to enumerating the number of graphs
locally equivalent to a fixed one. It will be shown that this
number is either
$\frac{(q-1)\epsilon(\mathcal{L})}{\lambda(\mathcal{L})}$ or its
half, depending on some issues discussed in the section.

Using this result, since we had estimated $\lambda(\mathcal{L})$,
the only remained step is to approximate $\epsilon(\mathcal{L})$,
the number of Eulerian vectors. This is done in section 7 using
{\it Tutte-Martin polynomial}. Indeed, $\epsilon(\mathcal{L})$,
number of Eulerian vectors can be written in term of this
polynomial. So, using the known recursive formula of Tutte-Martin
polynomial, we estimate $\epsilon(\mathcal{L})$.

In section 8, we put all these parts together to enumerate the
classes of local equivalency.

\section{Isotropic systems and graphic presentations}

Assume that $q$ is a power of $p$, which is an odd prime number,
and $\mathbf{F}_q$ is the finite field with $q$ elements. Also,
let $\mathbf{K}=\mathbf{F}^{2}_{q}$ denote a two-dimensional
vector space over $\mathbf{F}_q$, associated with the bilinear
form $\langle .,. \rangle$ satisfying $$\langle
(x,y),(x',y')\rangle=xy'-x'y$$ for every $(x,y),(x',y')\in
\mathbf{K}$. For a set $V$ of $n$ elements, define $\mathbf{K}^V$
to be the $2n$-dimensional vector space over $\mathbf{F}_q$,
equipped with the bilinear form $$\langle A,A'\rangle=\sum_{v\in
V} \langle A(v),A'(v)\rangle.$$

For $X, Y\in \mathbf{F}^V_q$, let $X\times Y$ be a vector in
$\mathbf{F}^V_q$ such that $(X\times Y)(v)=X(v)Y(v).$ Also, for
$v\in V$ and $(x,y)\in \mathbf{K}$, let $E_{v,(x,y)}$ be a vector
in $\mathbf{K}^V$ where its coordinates are all zero except the
$v$-th one which is equal to $(x,y)$, i.e.,
$E_{v,(x,y)}(w)=\delta_{vw}(x,y)$ for every $w\in V$. For any
vector $A\in \mathbf{K}^V$, let $A_v=E_{v,A(v)}$.

For simplicity, we present a vector $A$ in $\mathbf{K}^V$ as
$A=(X,Y)$, where $X, Y\in \mathbf{F}^V_q$. Therefore,
$A(v)=(X(v),Y(v))$. Also, for $X\in \mathbf{F}^V_q$ let $diag\ X$
be an $n\times n$ diagonal matrix where $(diag\ X)_{vv}=X(v)$,
and for the $2n\times 2n$ matrix
$$D=\begin{pmatrix}
  diag\ X & diag\ Y \\
  diag\ X' & diag\ Y'
\end{pmatrix},$$
define $$det\ D=diag X\ diag Y'-diag X'\ diag Y,$$ being an
$n\times n$ diagonal matrix. Therefore, $\langle A,
A'\rangle=tr(det\ D)$ where $tr$ is the usual trace function.
\footnote{Note that, the usual determinant of $D$ is equal to
$det(det\ D)$, which is the multiplication of diagonal entries of
$det\ D$.}

\begin{definition} An {\it isotropic system} $\mathcal{L}$ is an
$n$-dimensional subspace of $\mathbf{K}^V$ where $|V|=n$, such
that $\langle A, B\rangle=0$ for every $A, B\in \mathcal{L}$. In
other words, $\mathcal{L}$ is a subspace which is orthogonal to
itself.
\end{definition}

Note that, $\langle.,. \rangle$ is non-degenerate and since
$\mathcal{L}$ has dimension $n$, therefore $$\mathcal{L}=\{A\in
\mathbf{K}^V: \langle A, B \rangle=0,\ \forall B\in
\mathcal{L}\}.$$ In fact, if we fix a set of generators
$\{A_1=(X_1,Y_1),\dots A_n=(X_n,Y_n)\}$ for $\mathcal{L}$ and
construct the $n\times 2n$ matrix

\begin{equation}\label{basis}{\mathcal{B}=\begin{pmatrix}
   & A_{1} &  \\
   & A_{2} &  \\
   & \vdots &  \\
   & A_{n} &
\end{pmatrix}=\begin{pmatrix}
  X_{1} & \mid & Y_{1} \\
  X_{2} & \mid & Y_{2} \\
       & \vdots &  \\
  X_{n} & \mid & Y_{n}
\end{pmatrix},}\end{equation} then $\mathcal{L}$ is an isotropic system if and only if
$\mathcal{B}$ is a full-rank matrix (and hence $\dim
\mathcal{L}=n$) and
$$\mathcal{B}\begin{pmatrix}
  0 & I \\
  -I & 0
\end{pmatrix}\mathcal{B}^T=0,$$ Because, the $ij$-th entry of this
matrix is $$X_iY^T_j-Y_iX^T_j=\langle A_i,A_j\rangle.$$

We will shortly prove that every isotropic system can be defined
based on graphs. To start with, we introduce these graph-based
isotropic systems. Suppose that $G=(V,E)$ is a simple labeled
graph (without loops and multiple edges), where the label of each
edge comes from $\mathbf{F}_q$. Thus, we can represent the graph
by an $n\times n$ matrix $G=(g_{vw})$, where $n$ is equal to
$|V|$, the number of vertices in the graph, and for every $v,w\in
V $, $g_{vw}$ is equal to the label of the edge $vw$. So, $G$ is a
symmetric matrix with zero diagonal. Assume that $A=(X,Y)$ and
$B=(Z,T)$ are in $\mathbf{K}^V$ such that $diag Z\ diag Y-diag X\
diag T=cI $, where $I$ is the identity matrix and $c\in
\mathbf{F}_q $ is a non-zero constant. Denote by $\mathcal{L}$
the vector space generated by all vectors $g(v)(diag X\mid diag
Y)+B_v$ for $v\in V$, where $g(v)$ denotes the $v$-th row of $G$.
In fact, rows of matrix

\begin{equation}\label{igd}{(I \mid G)\cdot D }\end{equation}
form a basis for $\mathcal{L}$, where $$D=\begin{pmatrix}
  diag\ Z & diag\ T \\
  diag\ X & diag\ Y
\end{pmatrix}.$$

In order to prove that $\mathcal{L}$ is an isotropic system, first
note that $(I \mid G)$ is a full-rank matrix. Also, the
determinant of $D$ is $c^n$ which is non-zero, and hence $D$ is
full-rank as well, and since $D$ is a square matrix, $(I \mid
G)\cdot D$ is full rank too. On the other hand, we need to show
that the rows of $(I\mid G)\cdot D$ are orthogonal to each other,
or equivalently, $$(I\mid G)\cdot D \cdot \begin{pmatrix}
  0 & I \\
  -I & 0
\end{pmatrix}\cdot \big((I\mid G)\cdot D\big)^T=0.$$
But, the left hand side is equal to

\begin{eqnarray*}&=&(I\mid G)\begin{pmatrix}
  diag\ Z & diag\ T \\
  diag\ X & diag\ Y
\end{pmatrix}\begin{pmatrix}
  0 & I \\
  -I & 0
\end{pmatrix}\begin{pmatrix}
  diag\ Z & diag\ X \\
  diag\ T & diag\ Y
\end{pmatrix}\begin{pmatrix}
  I \\
  G
\end{pmatrix}\\
&=&(I\mid G)\begin{pmatrix}
  diag\ Z & diag\ T \\
  diag\ X & diag\ Y
\end{pmatrix}\begin{pmatrix}
  diag\ T & diag\ Y \\
  -diag\ Z & -diag\ X
\end{pmatrix}\begin{pmatrix}
  I \\
  G
\end{pmatrix}\\
&=&(I\mid G)\begin{pmatrix}
  0 & cI \\
  -cI & 0
\end{pmatrix}\begin{pmatrix}
  I \\
  G
\end{pmatrix}\\
&=&c(I\mid G)\begin{pmatrix}
  G \\
  -I
\end{pmatrix}=c(G-G)=0.
\end{eqnarray*}
Therefore, $\mathcal{L}$ is an isotropic system.

\begin{definition} Suppose that $\mathcal{L}$ is an $n$-dimensional
isotropic system for which there exist a graph $G$ and vectors
$A=(X,Y)$ and $B=(Z,T)$ in $K^V$ such that, $det\ D(A, B)=cI$ for
some $0\neq c\in \mathbf{F}_q$, where $$D(A,B)=
\begin{pmatrix}
  diag Z & diag T \\
  diag X & diag Y
\end{pmatrix},$$ and $(I\mid G)\cdot D(A,B)$ is a basis for $\mathcal{L}$. If such $G$ and $A, B$ exist,
we say $(G, A, B)$ is a {\it graphic presentation} of
$\mathcal{L}$, and $G$ is a {\it fundamental graph} of
$\mathcal{L}$.
\end{definition}

\begin{thm}\label{gpres}{ Every isotropic system $\mathcal{L}$ has a graphic
presentation.
 }\end{thm}

\begin{proof}{ We have already shown that every subspace which addmits a
graphic presentation is an isotropic system. Hence, it is
sufficient to prove that every isotropic system has a basis of
the form (\ref{igd}).

Consider an arbitrary basis for $\mathcal{L}$ and put them in the
rows of an $n\times 2n$ matrix $\mathcal{B}$ to get a matrix of
the form (\ref{basis}). Notice that, if we change the $v$-th
column of the first block of $\mathcal{B}$ with $v$-th column of
second block, we come up with a basis for another isotropic
system and it is equivalent to multiplying this matrix with a
$2n\times 2n$ matrix $D_1$ which consists of four $n\times n$
diagonal matrices (in fact, just two of these matrices are
non-zero). Among all such matrices $\mathcal{B}D_1$, choose the
one in which the rank of its first block is the maximum possible,
namely $r$. Now note that, in $K^V$ changing the order of
coordinates is equivalent to changing the order of columns of
$\mathcal{B}$ in the first and the second blocks. In fact, it is
equivalent to multiplying $\mathcal{B}$ by a $2n\times 2n$
permutation matrix from the right hand side, i.e., by a matrix of
the form

\begin{equation}\label{pi}{\Pi=\begin{pmatrix}
  \pi & 0 \\
  0 & \pi
\end{pmatrix},}\end{equation} where $\pi$ is a permutation matrix over $n$
elements. We find the permutation $\pi$ such that in
$\mathcal{B}D_1\Pi$ the first $r$ columns of the first block are
linearly independent. Then, there exists an invertible matrix $U$
such that
$$U\mathcal{B}D_1\Pi=\begin{pmatrix}
  I_r & \alpha & \mid & \beta & \gamma \\
  0 & \zeta & \mid & \eta & \theta
\end{pmatrix}.$$
Due to the properties of the matrix $D_1$ and the maximality
assumption, we have $\zeta=0$ and $\theta=0$. Since, rows of
matrix $U\mathcal{B}D_1\Pi$ form a basis for an isotropic system,
and because of the orthogonality assumption, we conclude that
$\eta=0$. Therefore, the rank of whole matrix is $r$. But this is
a basis for an isotropic system of dimension $n$. Thus, $r=n$, and
we have $$U\mathcal{B}=(I_n\mid
\beta)\Pi^{-1}D^{-1}_1=(\pi^{-1}\mid \beta\pi^{-1})D^{-1}_1.$$
Therefore, $\pi U\mathcal{B}=(I\mid \pi\beta\pi^{-1})D^{-1}_1$.
The matrix $\pi\beta\pi^{-1}$ may have non-zero diagonal entries,
but by multiplying $(I\mid \pi\beta\pi^{-1})$ by a $2n\times 2n$
matrix with four diagonal blocks, one can obtain a matrix with
the identity matrix in the first block, and the second block the
same as before, except that the diagonal entries are all zero.
Considering this multiplication, we end up with $\pi
U\mathcal{B}=(I\mid G')D'$, where $D'$ is the described $2n\times
2n$ matrix with four diagonal blocks. Note that both matrices
$\pi$ and $U$ are invertible, so that the rows of $\pi
U\mathcal{B}$ are still a basis for $\mathcal{L}$. Let
$$D'=\begin{pmatrix}
  diag Z' & diag T' \\
  diag X' & diag Y'
\end{pmatrix}.$$
By considering the orthogonality assumption, we conclude that
\begin{eqnarray*}0&=&(I\mid G')D'\begin{pmatrix}
  0 & I \\
  -I & 0
\end{pmatrix}D'^T(I\mid G')^T\\
&=&(I\mid G')\begin{pmatrix}
  diag Z' & diag T' \\
  diag X' & diag Y'
\end{pmatrix}\begin{pmatrix}
  0 & I \\
  -I & 0
\end{pmatrix}\begin{pmatrix}
  diag Z' & diag X' \\
  diag T' & diag Y'
\end{pmatrix}\begin{pmatrix}
  I \\
  G'^T
\end{pmatrix}\\
&=&(I\mid G')\begin{pmatrix}
  0 & det D' \\
  -det D' & 0
\end{pmatrix}\begin{pmatrix}
  I \\
  G'^T
\end{pmatrix}\\
&=&detD'\ G'^T-G'\ detD'.
\end{eqnarray*} Therefore, $detD'G'^T=G'\ detD'$. It means that the matrix
$G=G'\ detD'$ is symmetric. Moreover, it has a zero diagonal
since $G'$ does. We have

\begin{eqnarray*}\pi U\mathcal{B}&=&(I\mid
G)\begin{pmatrix}
  I & 0 \\
  0 & detD'^{-1}
\end{pmatrix}D'\\
&=&(I\mid G)\begin{pmatrix}
  diag Z' & diag T' \\
  detD'^{-1}diag X' & detD'^{-1}diag Y'
\end{pmatrix}.
\end{eqnarray*} Hence, if we define $$D=\begin{pmatrix}
  diag Z' & diag T' \\
  detD'^{-1}diag X' & detD'^{-1}diag Y'
\end{pmatrix},$$ then $\pi
U\mathcal{B}=(I\mid G)D$ and $det\ D=I$. This is a basis of the
form (\ref{igd}), and the proof is completed.

}\end{proof}

\section{Fundamental graphs of isotropic systems}

In the previous section we proved that every isotropic system
admits a fundamental graph. But, this fundamental graph is not
unique. In order to study these different fundamental graphs for
an isotropic system, we present a couple of definitions.

\begin{definition}  Let $G$ be a graph over the vertex set $V$. For
$v\in V$ and a number $r\in \mathbf{F}_q$, define $G*_{r} v$ to
be a graph (more precisely, a symmetric matrix with zero diagonal)
$G'=(g'_{uw})$, such that for every $w$, $g'_{vw}=g_{vw}$, and
also for every $u,w$ unequal to $v$,
$$g'_{uw}=g_{uw}+rg_{vu}g_{vw}.$$

Moreover, for a non-zero number $s\in \mathbf{F}_q$, define
$G\circ_s v$ to be a graph $G'=(g'_{uw})$, such that for each $u$,
$g'_{uv}=sg_{uv}$, and also for each $u,w$ unequal to $v$,
$g'_{uw}=g_{uw}$.

Two graphs $G$ and $G'$ are called {\it locally equivalent} if
there exists a sequence of the above operations, acting on $G$
obtains $G'$.
\end{definition}

Notice that, the operations $*$ and $\circ$ are invertible, so
that local equivalency is really an equivalency relation.

The following theorem is proved in \cite{eff}, and we do not
repeat the proof here.

\begin{thm}\label{loeq}{ Two graphs $G$ and $G'$ are locally
equivalent if and only if there exists an invertible $n\times n $
matrix $U$ and a $2n\times 2n$ matrix $D$ with four diagonal
blocks satisfying $det D=I$ and $U(I\mid G)D=(I\mid G')$.

\hfill{$\Box$}

 }\end{thm}

For a non-zero number $c\in \mathbf{F}_q$, let $G'=cG$ be the
usual product of a matrix $G$ and a constant $c$, i.e.,
$g'_{vw}=cg_{vw}$ for any $v, w$.

\begin{thm}\label{loeq1}{ For any two fundamental graphs $G$ and $G'$ of an
isotropic system, there exists a non-zero $c\in \mathbf{F}_q$ such
that $cG$ and $G'$ are locally equivalent. Conversely, if for
some non-zero number $c$, the graphs $cG$ and $G'$ are locally
equivalent, then there is an isotropic system such that $G$ and
$G'$ are its fundamental graphs.

}\end{thm}

\begin{proof}{ Suppose that $(G, A, B)$ and $(G', A', B')$ are two
graphic presentations for the isotropic system $\mathcal{L}$. It
means that, the rows of each of the matrices $(I\mid G)D(A, B)$
and $(I\mid G')D(A',B')$ form a basis for $\mathcal{L}$.
Therefore, there exists an invertible matrix $U$ such that
$$U(I\mid G)D(A, B)=(I\mid G')D(A',B').$$ Hence
$$U(I\mid G)D(A, B)D(A',B')^{-1}=(I\mid G').$$

Note that, both $det D(A, B)$ and $det D(A', B')$ are (non-zero)
constant numbers, and $det D(A', B')^{-1}$ is also a constant.
Thus, there exists a non-zero number $c\in \mathbf{F}_q$ such that
$det \big(D(A, B)D(A',B')^{-1}\big)=cI$. Now let
$$D=\begin{pmatrix}
  I & 0 \\
  0 & c^{-1}I
\end{pmatrix}D(A, B)D(A',B')^{-1}.$$ $det D=I$ and we have $$U(I\mid G)\begin{pmatrix}
  I & 0 \\
  0 & cI
\end{pmatrix}D=(I\mid G').$$ Then, $U(I\mid cG)D=(I\mid G')$
and the first part of the conclusion follows from theorem
\ref{loeq}.

Conversely, suppose that $cG$ and $G'$ are locally equivalent.
Therefore, there exist matrices $U$ and $D$ such that $U$ is
invertible, $det D=I$ and $U(I\mid cG)D=(I\mid G')$. Also,
suppose that $A$ and $B$ are two vectors such that $D=D(A,B)$.
Therefore, $(G,cA, B)$ and $(G', A', B')$ are two graphic
presentations for the same isotropic system, where $A',B'$ are
defined so that $D(A', B')=I_{2n}$. More precisely, $A'(v)=(0,1)$
and $B'(v)=(1,0)$ for each $v$.}

\end{proof}

Having this theorem in hand, we can now study the classes of
local equivalency of graphs by investigating different graphic
presentations of an isotropic system.

\section{Eulerian vectors and local complementation}

We call a vector $A\in \mathbf{K}^V$ {\it complete} if $A(v)$ is
non-zero for all $v\in V$.  Assume that $A\in \mathbf{K}^V$ is
complete. We define $\hat{A}$ to be the vector subspace generated
by vectors $A_v$ for all $v\in V$, i.e., $\hat{A}=\langle A_v:
v\in V\rangle$. In fact, if $A=(X,Y)$ then $\hat{A}$ is the vector
space generated by rows of $(diag X, diag Y)$.

\begin{definition} Let $\mathcal{L}$ be an isotropic system and $A\in
\mathbf{K}^V$ be a complete vector. We call $A$ an {\it Eulerian
vector} for $\mathcal{L}$ if $\hat{A}\cap \mathcal{L}=0$.
\end{definition}

\begin{lem}\label{evgr}{ Suppose that $\mathcal{L}$ is an isotropic system and $(G, A,
B)$ is a graphic presentation for $\mathcal{L}$. Then $A$ is an
Eulerian vector for $\mathcal{L}$.

}\end{lem}

\begin{proof}{ By definition of the graphic presentation we know that
$det D(A, B)$ is a non-zero constant. Therefore, $A(v)$ is
non-zero for any $v$. In fact, the $v$-th entry of the diagonal of
$det D(A, B)$ is $\langle B(v), A(v)\rangle$, so that $A(v)$
should be non-zero and then $A$ is complete. To get a
contradiction, suppose that $\hat{A}\cap \mathcal{L}$ is
non-zero. Therefore, if $A=(X, Y)$ then there exists a non-zero
vector $S\in \mathbf{F}_q^{n}$ such that $S(diag X | diagY)\in
\mathcal{L}$. $S$ is non-zero, hence at least one of its
coordinates , say $v$-th, is non-zero. By considering the
orthogonality condition, we have

\begin{eqnarray*} 0 &=&\langle S(diagX | diag
Y), g(v)(diag X | diag Y)+B_v\rangle \\
&=&\sum_{w\in V} S(w)g_{vw} \langle(X(w), Y(w)),(X(w),
Y(w))\rangle+\langle S(diag X | diag Y), B_v \rangle\\
& = &S(v)\langle A(v),B(v)\rangle.
\end{eqnarray*}
But, $S(v)$ is non-zero and also $\langle A(v),B(v)\rangle$ is
non-zero since $det D(A, B)$ is so, which is a contradiction.

}\end{proof}

\begin{cor}\label{eise}{ Every isotropic system has an Eulerian
vector.

}\end{cor}

\begin{proof}{ By theorem \ref{gpres}, every isotropic system
has a graphic presentation $(G, A, B)$ and by lemma \ref{evgr},
$A$ is an Eulerian vector for the isotropic system.

}\end{proof}

Let $\mathcal{L}$ be an isotropic system and $A$ an Eulerian
vector for $\mathcal{L}$. Therefore $\mathcal{L}\cap \hat{A}=0$.
If $A'$ is a vector which is equal to $A$ at any coordinate
except the $v$-th one, at which it is equal to a non-zero multiple
of $A(v)$, then $\hat{A'}=\hat{A}$, and therefore $A'$ is also an
Eulerian vector for $\mathcal{L}$.

This observation gives us the motivation of defining
$\mathbf{K}^*$ to be $\mathbf{K}\backslash\{0\}$ under the
equivalency relation $(x,y)\sim(x',y')$ iff $(x,y)=r(x',y')$, for
some non-zero $r$, (and hence $\mid \mathbf{K}^*\mid =q+1$). Now
by the above discussion if we replace each coordinate of $A$ with
something equivalent to it, we obtain another Eulerian vector. The
set of Eulerian vectors of an isotropic system has even more
useful properties.

\subsection{Switching property}

\begin{definition} We say that a subset $\sum\subseteq \mathbf{K}^V$ of
complete vectors has the {\it switching property} if

\vspace{3mm}

\noindent{\rm (i)} Similar to Eulerian vectors, for each $A\in
\sum$, one can replace each coordinate of $A$ by its (scalar)
multiple and it still remains in this subset.

\vspace{3mm}

\noindent{\rm (ii)} In addition, for each $A\in \sum$ and $v\in
V$, $A-A_v+rE_{v,(x,y)}$ is still in $\sum$ for each non-zero
$r\in \mathbf{F}_q$ and every $(x,y)\in \mathbf{K}^*$ except one
$(x,y)$. In other words, in a set with switching property and a
vector $A$ in this set, we can replace a coordinate of $A$ with
exactly $q$ elements of $\mathbf{K}^*$ so that it still remains in
$\sum$.
\end{definition}

We will observe shortly that switching at the vertex $v$ is
equivalent to a local complementation operation on this vertex.

\begin{thm}\label{switch}{ The set of Eulerian vectors of an isotropic
system has the switching property.

}\end{thm}

\begin{proof}{ Assume that $A$ is an Eulerian vector for
the isotropic system $\mathcal{L}$, and to lead to a
contradiction, suppose that $(x_i,y_i)$, $i=1,2$, are two
different vectors in $\mathbf{K}^*$ and $A_i=A-A_v+E_{v,(x_i,
y_i)}$, $i=1,2$, are not Eulerian (we know that all of these
vectors $A-A_v+E_{v,(x, y)}, (x,y) \in \mathbf{K}^* $ can not be
Eulerian, since if so, for every $C \in \mathcal{L}, C(v)=0$ and
the dimension of $\mathcal{L}$ could not be equal to $n$).
Therefore there exist non-zero vectors $C_i\in \hat{A_i}\cap
\mathcal{L}$ for $i=1,2$. The $v$-th coordinate of $C_i$ can not
be zero, because otherwise, $C_i\in \hat{A}\cap \mathcal{L}$,
which is not possible since $A$ is Eulerian. Therefore, the $v$-th
coordinate of $C_i$ is a non-zero multiple of $(x_i,y_i)$,
$i=1,2$. Now notice that $(x_1,y_1)$ and $(x_2,y_2)$ are
different elements of $\mathbf{K}^*$, thus they are linearly
independent and there exist $r_1,r_2\in \mathbf{F}_q$ such that
$A(v)=r_1(x_1,y_1)+r_2(x_2,y_2)$. Therefore $r_1C_1+r_2C_2$ is a
non-zero vector in $\hat{A}\cap \mathcal{L}$, which is a
contradiction.

}\end{proof}

\begin{thm}\label{evgp}{ For every isotropic system $\mathcal{L}$ and a
graphic presentation $(G, A, B)$ of it, $A$ is an Eulerian
vector. Conversely, for every Eulerian vector $A$, there exists a
graphic presentation $(G, A, B)$. Also, this graphic presentation
is unique up to a (non-zero) constant, i.e., if $(G', A, B')$ is
another graphic presentation for $\mathcal{L}$ then there exists
a non-zero number $c\in \mathbf{F}_q$ such that $G'=cG$ and
$B'=cB$. }\end{thm}

\begin{proof}{The first part of the theorem was already proved in lemma
\ref{evgr}. For the second part, suppose that $A=(X,Y)$ is an
Eulerian vector of the isotropic system $\mathcal{L}$. By the
switching property, for every $v\in V$, there exists some $(z_v,
t_v)\in \mathbf{K}$ such that $(z_v, t_v)\nsim A(v)$ (meaning
that $(z_v, t_v)$ and $A(v)$ are not scalar multiples of each
other) and there is a vector of the form $C_v+E_{v,(z_v,t_v)}$ in
$\mathcal{L}$, where $C_v$ is in $\hat{A}$ and $C_v(v)=0$. Since
$(z_v, t_v)\nsim A(v)$, we have $\langle (z_v,t_v), A(v)
\rangle\neq 0$ and by considering an scalar multiple (if
necessary) of $(z_v,t_v)$, we may assume that

\begin{equation}\label{d1}{\langle (z_v,t_v), A(v)
\rangle=1.}\end{equation}

We know that $C_v\in \hat{A}$ and $C_v(v)=0$, thus there exists a
matrix $G=(g_{uw})$ over $\mathbf{F}_q$ such that $g_{vv}=0$ and
$C_v=g(v)(diagX | diag Y)$, for every $v$. Here, by $g(v)$ we mean
the $v$-th row of $G$.

Now, Let $B=(Z, T)$ be the vector in $\mathbf{K}^V$ with
$Z(v)=z_v$ and $T(v)=t_v$ for every $v\in V$. Due to the equation
(\ref{d1}), we observe that $det D(A, B)=I$ and also $D(A, B)$ is
an invertible matrix. Using this notation, the rows of $(I\mid
G)D(A, B)$ are all in $\mathcal{L}$. Since, this matrix is a
full-rank one, its rows form a basis for $\mathcal{L}$.
Therefore, once we show that $G$ is the matrix for a graph, we
will end up with the presentation $(G, A, B)$ for $\mathcal{L}$.

To show that $G$ is a graph, first note that $g_{vv}=0$ by its
definition. For proving that $G$ is symmetric, consider again the
orthogonality assumption. We have

\begin{eqnarray*} 0&=&(I\mid G)\begin{pmatrix}
  diag Z & diag T \\
  diag X & diag Y
\end{pmatrix}\begin{pmatrix}
  0 & I \\
  -I & 0
\end{pmatrix}\begin{pmatrix}
  diag Z & diag X \\
  diag T & diag Y
\end{pmatrix}\begin{pmatrix}
  I \\
  G^T
\end{pmatrix}\\
&= &(I\mid G)\begin{pmatrix}
  diag Z & diag T \\
  diag X & diag Y
\end{pmatrix}\begin{pmatrix}
  diag T & diag Y \\
  -diag Z & -diag X
\end{pmatrix}\begin{pmatrix}
  I \\
  G^T
\end{pmatrix}\\
&=&(I\mid G)\begin{pmatrix}
  0 & I \\
  -I & 0
\end{pmatrix}\begin{pmatrix}
  I \\
  G^T
\end{pmatrix}\\
&=&G^T-G.
\end{eqnarray*}
Therefore, $G$ is a graph, and $(G,A, B)$ is a graphic
presentation of $\mathcal{L}$.

For the uniqueness, suppose that $(G', A, B')$ is another graphic
presentation for $\mathcal{L}$ such that $det D(A, B')=cI$. It
means that $(I\mid G')D(A, B')$ is also a basis for
$\mathcal{L}$. Let $B'=(Z', T')$ and by considering the
orthogonality assumption, once again we have

\begin{eqnarray*}0&=&(I\mid G)\begin{pmatrix}
  diag Z & diag T \\
  diag X & diag Y
\end{pmatrix}\begin{pmatrix}
  0 & I \\
  -I & 0
\end{pmatrix}\begin{pmatrix}
  diag Z' & diag X \\
  diag T' & diag Y
\end{pmatrix}\begin{pmatrix}
  I \\
  G'
\end{pmatrix}\\
&=&(I\mid G)\begin{pmatrix}
  det D(B', B) & det D(A, B) \\
  -det D(A, B') & 0
\end{pmatrix}\begin{pmatrix}
  I \\
  G'
\end{pmatrix}\\
&=&(I\mid G)\begin{pmatrix}
  det D(B', B)+G' \\
  -cI
\end{pmatrix}\\
&=& det D(B', B)+G'-cG.
\end{eqnarray*}

Therefore, $det D(B', B)+G'-cG=0$ and since $det D(B', B)$ is a
diagonal matrix and the diagonals of $G$ and $G'$ are both equal
to zero, we have $det D(B', B)=0$ and $G'=cG$. Hence, it just
remains to show $B'=cB$. Because of the equation $$0=det D(B', B)=
diag Z\ diag T'-diag Z'\ diag T,$$ for any $v$ there exists a
$c_v\in \mathbf{F}_q$ such that $(Z'(v), T'(v))=c_v(Z(v), T(v))$.
Therefore $D(A, B')_{vv}=c_vD(A, B)_{vv}$. On the other hand $det
D(A, B)=I$ and $det D(A, B')=cI$, so that $c_v=c$ for any $v$.
Hence $B'=cB$.

}\end{proof}

Using this theorem, if we denote by $\epsilon(\mathcal{L})$ the
number of Eulerian vectors of the isotropic system $\mathcal{L}$,
we conclude the following corollary.

\begin{cor}\label{eulpres}{ The number of graphic presentations of an
isotropic system is equal to $(q-1)\epsilon(\mathcal{L})$.

}\end{cor}

\subsection{Switching property in terms of local complementation}

The following theorem, explains the relationship between the
switching property and local complementation.

\begin{thm}\label{locswit}{ Suppose that $(G, A, B)$ is a graphic
presentation for the isotropic system $\mathcal{L}$, and $v\in V$.

\vspace{3mm}

\noindent{\rm (i)} If $r \in \mathbf{F}_q$, then $(G*_rv,
A+rB_v,B+rg(v)^2\times A)$ is also a graphic presentation of
$\mathcal{L}$. Therefore switching $A$ at $v$ is equivalent to a
local complementation operator.

\vspace{3mm}

\noindent{\rm (ii)} If $s\in \mathbf{F}_q$ is non-zero, then
$(G\circ_sv, A+(s^{-1}-1)A_v, B+(s-1)B_v )$ is also a graphic
presentation of $\mathcal{L}$.

}\end{thm}

\begin{proof}{ We prove first part, and the second part is similar. It is easy to
check that $det D(A+rB_v,B+rg(v)^2\times A)$ is constant. Hence,
it is sufficient to show that all of the rows of $(I\mid G*_rv)D(
A+rB_v,B+rg(v)^2\times A)$ are in $\mathcal{L}$. Let $G'=G*_rv$
and $w\in V$, $w\neq v$. We have
$g'(w)=g(w)+rg_{vw}g(v)-rg^2_{vw}\delta_w$, thus the $w$-th row of
$(I\mid G')D(A+rB_v,B+rg(v)^2\times A)$ is equal to

\begin{eqnarray*}
&= & g'(w)\times (A+rB_v)+(B_w+rg^2_{vw}A_w)\\
&= &(g(w)+rg_{vw}g(v)-rg^2_{vw}\delta_w)\times
(A+rB_v)+(B_w+rg^2_{vw}A_w)\\
&= & g(w)\times A+rg_{vw}g(v)\times
A-rg^2_{vw}A_w+rg_{vw}B_v+B_w+rg^2_{vw}A_w \\
&= &\big(g(w)\times A+B_w\big)+rg_{vw}\big(g(v)\times
A+B_v\big),\\
\end{eqnarray*}
which is in $\mathcal{L}$. Also, for the $v$-th row, $g'(v)=g(v)$
and
$$g(v)\times (A+rB_v)+(B_v+rg^2(v)\times A_v)= g(v)\times A+B_v$$
which is an element of $\mathcal{L}$.

}\end{proof}

\section{Index of an isotropic system}

We are now in the position of introducing the notion of {\it
index} for an isotropic system.

\begin{thm}\label{index1}{ For any isotropic system $\mathcal{L}$, there
exists a number $\lambda(\mathcal{L})$ such that for any
fundamental graph $G$ for $\mathcal{L}$, there are exactly
$\lambda(\mathcal{L})$ pairs $(A, B)$ such that $(G, A, B )$ is a
graphic presentation of $\mathcal{L}$. This number is called the
index of the isotropic system $\mathcal{L}$.
 }\end{thm}

\begin{proof}{ Suppose that $G$ is a fundamental graph of $\mathcal{L}$, and
there exist exactly $k$ graphic presentations of the form $(G,
A_i, B_i)$, $i=1, \dots ,k$ for $\mathcal{L}$. It is sufficient
to show that for any other fundamental graph $H$ for
$\mathcal{L}$, there are also $k$ graphic presentation with $H$
as a fundamental graph. Since, $G$ and $H$ are fundamental graphs
of the same isotropic system, by theorems \ref{loeq} and
\ref{loeq1}, there exist invertible matrices $U$ and $D$, such
that $D$ consists of four diagonal blocks and $det D=cI$ for some
non-zero $c\in \mathbf{F}_q$, and moreover,
\begin{equation}\label{l2}{U(I \mid H)D=(I\mid G).}\end{equation}

Since, $(G, A_i, B_i)$ is a graphic presentation of $\mathcal{L}$
for $i=1, \dots, k$, the rows of $(I\mid G)D(A_i, B_i)$ form a
basis for $\mathcal{L}$. Now using (\ref{l2}), we conclude that
the rows of $U(I\mid H)D\ D(A_i, B_i)$ and hence the rows of
$(I\mid H)D\ D(A_i, B_i)$ form bases for $\mathcal{L}$.

Notice that $det D\ D(A_i, B_i)$ is a constant. Therefore, it
gives us a graphic presentation of $\mathcal{L}$ with fundamental
graph $H$, for $i=1, \dots ,k$. Also, since the matrices $U, D$
and $D(A_i, B_i)$ are invertible, the described $k$ presentations
are different. Moreover, for any presentation with fundamental
graph $H$, we can convert it to a presentation with fundamental
graph $G$. Thus, for any fundamental graph of $\mathcal{L}$, there
exist exactly $\lambda(\mathcal{L})=k$ graphic presentations with
this graph as a fundamental graph.

}\end{proof}

\begin{thm}\label{index2}{ Assume that $\mathcal{L}$ is an isotropic system admitting $G$ as a
fundamental graph. Then $\lambda(\mathcal{L})$ is equal to the
number of matrices of the form $D(A, B)$, with non-zero constant
determinant, and
\begin{equation}\label{intsol}{(I\mid G)D(A, B)\begin{pmatrix}
  G \\
  -I
\end{pmatrix}=0.}\end{equation} In fact, $\lambda(\mathcal{L})=\lambda(G)$, meaning that the index of an isotropic
system just depends on any arbitrary fundamental graph.

}\end{thm}

\begin{proof}{ As in the proof of theorem \ref{index1}, suppose that $(G, A_i,
B_i)$, $i=1, \dots, k$, are all graphic presentations of
$\mathcal{L}$ with fundamental graph $G$. Using the orthogonality
assumption we have

\begin{eqnarray*}0 & = & (I\mid G)D(A_1, B_1)\begin{pmatrix}
  0 & I \\
  -I & 0
\end{pmatrix}D(A_i, B_i)^T\begin{pmatrix}
  I \\
  G
\end{pmatrix}\\
&=&(I\mid G)D(A_1, B_1)\begin{pmatrix}
  0 & I \\
  -I & 0
\end{pmatrix}D(A_i, B_i)^T \begin{pmatrix}
  0 & -I \\
  I & 0
\end{pmatrix}\begin{pmatrix}
  G \\
  -I
\end{pmatrix}.
\end{eqnarray*}
On the other hand, the matrix $$\begin{pmatrix}
  0 & I \\
  -I & 0
\end{pmatrix}D(A_i, B_i)^T \begin{pmatrix}
  0 & -I \\
  I & 0
\end{pmatrix}$$ consists of four diagonal matrices and has
non-zero constant determinant, and also, $$D(A_1,
B_1)\begin{pmatrix}
  0 & I \\
  -I & 0
\end{pmatrix}D(A_i, B_i)^T \begin{pmatrix}
  0 & -I \\
  I & 0
\end{pmatrix}$$ satisfies all of these properties. Therefore, for each $i$, where
$i=1, \dots, k$, we find a solution of (\ref{intsol}).
Conversely, since all of the above equations can be inverted, for
each solution $D(A, B)$ we can find one of the graphic
presentations $(G, A_i, B_i)$.

}\end{proof}

\begin{cor}\label{indexpres}{ For an isotropic system that admits
$G$ as a fundamental graph, the number of its graphic
presentations is equal to $\lambda(G)$ times the number of graphs
that are locally equivalent to $cG$ for some non-zero $c\in
\mathbf{F}_q$.

}\end{cor}

%In chapter 3, the solutions of (\ref{intsol}) with determinant
%$I$ are called {\it internal solutions}, and some significant
%properties are proved for these solutions.

\subsection{Bineighborhood space and index of a graph}

For a graph $G$, we call a pair $vw$ of vertices an {\it edge}, if
$g_{vw}\neq 0$. Suppose that, $C$ is an even cycle (a cycle with
an even number of edges) consisting of vertices $v_1, v_2, \dots
,v_{2l}$. Set
$$\nu (C)=\sum^{2l}_{i=1} (-1)^{i}g_{v_i v_{i+1}} g(v_i)\times
g(v_{i+1}).$$

\begin{definition}  Suppose that $G$ is a graph. The {\it bineighborhood
space} of $G$, denoted by $\nu (G)$, is a subspace of
$\mathbf{F}^V_q$ defined by
$$\nu (G)= span  \{\nu(C): C\ even\ cycle\}\cup \{g(v)\times g(w):
g_{vw}=0\} .$$
\end{definition}

We assume that the graphs we consider are connected, and restate a
couple of theorems, which will be used shortly. These theorems
are all proved in \cite{eff}.

\begin{lem}\label{iden}{ If $D(A, B)$, where $A=(X, Y)$ and $B=(Z,
T)$, satisfies (\ref{intsol}) then $Y+Z$ is a scalar multiple of
$(1,1,\dots, 1)$. On the other hand, for any such vectors $A, B$,
the matrix $D(A, B)+cI_{2n}$ satisfies (\ref{intsol}) for each
$c\in \mathbf{F}_q$.

\hfill{$\Box$} }\end{lem}

\begin{thm}\label{binsp}{ Suppose that $D(A, B)$ satisfies
(\ref{intsol}) and  $A=(X, Y)$. Then $X\in \nu (G)^{\perp}$. On
the other hand for any $X\in \nu (G)^\perp$,

\vspace{3mm}

\noindent{\rm (i)} if $G$ has an odd cycle, then there exist a
unique $Y$ and a unique $T$ such that $D(A, B)$ satisfies
(\ref{intsol}), where $A=(X, Y)$ and $B=(-Y, T)$.

\vspace{3mm}

\noindent{\rm (ii)} if $G$ does not have an odd cycle, then there
exist exactly $q$ pairs of $Y_i, T_i$, $i=1,\dots, q$, such that
$D(A_i, B_i)$ satisfies (\ref{intsol}), where $A=(X, Y_i)$ and
$B=(-Y_i , T_i)$.

\hfill{$\Box$} }\end{thm}

\begin{thm}\label{cdet}{ If $D(A,B)$ satisfies (\ref{intsol}) for some vectors $A$
and $B$, then it has a constant determinant.

\hfill{$\Box$}}\end{thm}

Having all of these theorems in hand, we conclude the following
statement.

\begin{thm}\label{cont}{ \hfill

\vspace{3mm}

\noindent{\rm (i)} If $G$ has an odd cycle, then for any $X\in \nu
(G)^\perp$ there are exactly $q$ ($A, B$)'s, where the first
component of $A$ is $X$, and $D(A, B)$ satisfies (\ref{intsol}).
Among all of these $q$ solutions, either all of them or $q-2$ of
them have non-zero constant determinant.

\vspace{3mm}

\noindent {\rm (ii)} If $G$ has no odd cycle, then for any $X\in
\nu (G)^\perp$ there are exactly $q^2$ ($A, B$)'s where the first
component of $A$ is $X$ and $D(A, B)$ satisfies (\ref{intsol}).
Among all of these $q^2$ solutions,  at least $q(q-2)$ of them
have non-zero constant determinant.

}\end{thm}

\begin{proof}{ Suppose that $G$ has an odd cycle, and fix some $X\in \nu
(G)^\perp$. By lemma \ref{iden} and theorem \ref{binsp}, there are
$A_0=(X, Y_0)$ and $B_0=(-Y_0, T)$ such that $D(A_0, B_0)$
satisfies (\ref{intsol}), and any other solution $D(A, B)$ of
(\ref{intsol}), where the first component of $A$ is $X$, is of
the form $D(A, B)=D(A_0, B_0)+cI_{2n}$, for any $c\in
\mathbf{F}_q$. Hence, there are $q$ solutions with this property.
Notice that, $det D(A, B)= det D(A_0, B_0)+ c^2I$, and by theorem
\ref{cdet}, $\ det D(A_0, B_0)=d_0I$ is constant. Thus, depending
on whether $-d_0\in \mathbf{F}_q$ is a perfect square or not,
there are either $q-2$ or $q$ different values of $c\in
\mathbf{F}_q$ such that $d_0+c^2$ is non-zero.

The proof of {\rm (ii}) is exactly the same. The only difference
is that in this case there are $q$ solutions of the form $A=(X,
Y),\ B=(-Y, T)$ for (\ref{intsol}).

}\end{proof}

We can now give an estimation on $\lambda(G)$ for a graph $G$.

\begin{cor}\label{nsol}{ If a graph $G$ has an odd cycle
then $$(q-2)q^{dim\,\nu (G)^\perp}\leq \lambda(G)\leq q^{dim\,\nu
(G)^\perp+1 },$$ and if not $$(q-2)q^{dim\,\nu (G)^\perp +1}\leq
\lambda(G)\leq q^{dim\,\nu (G)^\perp+2 }.$$

}\end{cor}

\section{ The number of graphs locally equivalent to a fixed one}

In this section we plan to compute $l(G)$, the number of graphs
which are locally equivalent to the graph $G$. Once again, we
assume that all the graphs we consider are connected.

\begin{lem}\label{badih}{ $l(G)=l(cG)$ for any non-zero $c\in \mathbf{F}_q$.

}\end{lem}

\begin{proof}{ It is not hard to see that, for any graph $H$ locally equivalent to
$G$, $cH$ is locally equivalent to $cG$.

 }\end{proof}

\begin{lem}\label{lcl}{ If $c\in \mathbf{F}_q$ is a non-zero perfect square,
 then $G$ and $cG$ are locally equivalent.

}\end{lem}

\begin{proof}{ Suppose that $c=d^2$, and apply the operation $\circ_d v$
on $G$ for all $v\in V$, and the resulted graph is $cG$.

}\end{proof}

\begin{thm}\label{lcl2} { The number of graphs locally
equivalent to some non-zero scalar multiple of the graph $G$ is
either equal to $l(G)$ or $2l(G)$.

}\end{thm}

\begin{proof}{ Suppse that $cG$ is locally equivalent to $G$ for any non-zero
$c$. Then for any graph locally equivalent to $cG$ for some $c$,
it is also locally equivalent to $G$. Therefore $l(G)$ is the
number of graphs that are locally equivalent to $cG$ for some $c$.

Next, suppose that there exists a non-zero $c_0$ such that $c_0G$
is not locally equivalent to $G$. By lemma \ref{lcl}, $c_0$ is not
a perfect square. Assume that, $H$ is a graph that is locally
equivalent to $cG$, for some $c$. If $c$ is a perfect square,
then by lemma \ref{lcl}, $H$ is also locally equivalent to $G$.
On the other hand, if $c$ is not a perfect square, then
${c_0}^{-1} c$ is a perfect square. Hence ${c_0}^{-1}H$ is
locally equivalent to $G$. Therefore, for any such $H$, it is
locally equivalent to either $G$ or $c_0G$. Also, by lemma
\ref{badih}, $l(c_0G)=l(G)$. Therefore, in this case, the number
of graphs locally equivalent to $cG$ for some $c$, is $2l(G)$.

 }\end{proof}

Using theorem \ref{lcl2}, we can relate $l(G)$ to the number of
graphs that are locally equivalent to $cG$ for some non-zero $c$.
Since, this number appears in counting the number of graphic
presentations of an isotropic system, we can obtain some useful
information about $l(G)$. For this purpose, fix an isotropic
system $\mathcal{L}$ admitting $G$ as a fundamental graph.

\begin{thm}\label{kerel}{ The number of graphs locally equivalent
to $cG$ for some non-zero $c\in \mathbf{F}_q$ is equal to
$$\frac{(q-1)\epsilon(\mathcal{L})}{\lambda(\mathcal{L})}.$$

}\end{thm}

\begin{proof}{ By corollary \ref{eulpres}, the number of graphic
presentations of $\mathcal{L}$ is equal to
$(q-1)\epsilon(\mathcal{L})$. On the other hand, by corollary
\ref{indexpres}, the number of graphic presentations is equal to
$\lambda(\mathcal{L})$ times the number of graphs that are locally
equivalent to $cG$, for some non-zero $c$. By letting these two
values be equal, one obtains the described conclusion.

}\end{proof}

The following corollary is a direct consequence of theorems
\ref{lcl2} and \ref{kerel}.

\begin{cor}\label{kerel2}{ Either $l(G)$ or $2l(G)$ is equal to
$$\frac{(q-1)\epsilon(\mathcal{L})}{\lambda(\mathcal{L})}.$$

}\end{cor}

\vspace{4mm}

Using corollaries \ref{nsol} and \ref{kerel2}, giving a bound for
$\epsilon(\mathcal{L})$ can lead us to a bound for $l(G)$.

\vspace{3mm}

\noindent{\bf Remark.} {\it In \cite{bouchalg},\cite{bouchtree}
and \cite{bouchet}, it has been shown that in the binary case,
$l(G)=\frac{\epsilon(\mathcal{L})}{\lambda(\mathcal{L})}.$ So,
corollary \ref{kerel2} is valid for binary case, too.}

\section{The number of Eulerian vectors}

In this section, we include the binary case. To be precise, we
assume that $q$ is either $2$ or a power of an odd prime number.
In addition, we do {\it not} assume the connectivity of graphs. In
order to compute $\epsilon(\mathcal{L})$, the number of Eulerian
vectors of an isotropic system, we use a well-known polynomial,
so called {\it Tutte-Martin polynomial} which is defined for an
isotropic system $\mathcal{L}$ as follows:

$$ \mathcal{M}(\mathcal{L}; t)=\sum_{C:\ complete}  (t-q)^{\ dim(\mathcal{L}\cap \hat{C})},$$ where the summation is
over all complete vectors $C\in \mathbf{K}^V$. By the definition
of Eulerian vector,
$$\epsilon(\mathcal{L})=\mathcal{M}(\mathcal{L};\,q).$$ So, in order to
compute $\epsilon(\mathcal{L})$, it suffices to be able to
compute the Tutte-Martin polynomial. To do that, we give a
recursion formula for this polynomial.

\begin{definition} Let $\mathcal{L}$ be an isotropic system, $v\in V$ and
$x\in \mathbf{K}^*$ (or $x=0$). Define
$$\mathcal{L}^v_x=\{C\in \mathcal{L}: \langle C(v), x
\rangle=0\}=\{C\in \mathcal{L}: C(v)\sim x\ or\ C(v)=0 \},$$ and

$$\mathcal{L}\vert^v_x= projection\ of\ \mathcal{L}^v_x\ on\ \mathbf{K}^{V-\{v\}}.$$
Also, let $\mathcal{L}^v_0=\{C\in \mathcal{L}: C(v)=0\}$, and for
$C\in \mathbf{K}^V$, we set $C\vert^v$ to be the projection  of
$C$ on $\mathbf{K}^{V-\{v\}}$.
\end{definition}

\begin{lem}\label{minor}{ $\mathcal{L}\vert^v_x$ is an isotropic system.

}\end{lem}

\begin{proof}{ First notice that, $\mathcal{L}^v_x$ is a subspace of $\mathcal{L}$ and is
self-orthogonal. On the other hand, the $v$-th coordinate of each
vector in $\mathcal{L}^v_x$ is either $0$ or equivalent to $x$.
Then the $v$-th coordinate of any two vectors in $\mathcal{L}^v_x$
are orthogonal, and hence if we delete the $v$-th coordinate, all
of the vectors remain orthogonal to each other. In other words,
$\mathcal{L}\vert^v_x$ is again self-orthogonal. Thus, it remains
to show $dim\,\mathcal{L}\vert^v_x=n-1$. We consider two cases:

\vspace{3mm}

\noindent{\rm (i)} There exists $C_0\in \mathcal{L}$ such that
$\langle C_0(v), x \rangle\neq 0$. In this case, we have
$\mathcal{L}^v_x\neq \mathcal{L}$, and therefore, $dim
\mathcal{L}^v_x=n-1$. Notice that, $E_{v,x}$ is not orthogonal to
$C_0$, therefore $E_{v,x}$ is not in $\mathcal{L}^v_x$. Then, the
projection of $\mathcal{L}^v_x$ on $\mathbf{K}^{V-\{v\}}$ is
injective, and $dim \mathcal{L}\vert^v_x$ is also $n-1$.

\vspace{3mm}

\noindent{\rm (ii)} For any $C\in \mathcal{L}$, $\langle C(v),
x\rangle=0$. Then $E_{v,x}$ is in $\mathcal{L}$, and for any $C\in
\mathcal{L}$, $C(v)=0$ or $C(v) \sim x$. Thus,
$\mathcal{L}^v_x=\mathcal{L}$, and
$\mathcal{L}=\mathcal{L}^v_0\oplus \langle E_{v,x}\rangle$. It
means that, by removing the $v$-th coordinate of
$\mathcal{L}^v_0$, we end up with $\mathcal{L}\vert^v_x$.
Therefore, $dim \mathcal{L}\vert^v_x=dim \mathcal{L}^v_0=n-1$.

}\end{proof}

\vspace{3mm}

\noindent{\bf Remark.} {\it Let $G$ be a fundamental graph of
$\mathcal{L}$ and $(G, A, B)$ be a graphic presentation. Then,
$\mathcal{M}(\mathcal{L};t)$ just depends on $G$. In fact, if
$\mathcal{L}'$ is another isotropic system with a graphic
presentation $(G,A', B')$, then multiplication by the matrix
$D(A, B)^{-1}D(A', B')$ maps $\mathcal{L}$ to $\mathcal{L}'$.
Also, it maps $\mathcal{L}\cap \hat{C}$ to $\mathcal{L}'\cap
\hat{C'}$, where $C'=CD(A, B)^{-1}D(A', B')$. Therefore,
$\mathcal{M}(\mathcal{L'};t)=\mathcal{M}(\mathcal{L'};t)$.}

\vspace{3mm}

Using this remark, we can talk about the Tutte-Martin polynomial
of a graph, $\mathcal{M}(G;t)$. Also, if $G$ and $H$ are two
locally equivalent graphs, they are fundamental graphs of the
same isotropic system, and then
$\mathcal{M}(G;t)=\mathcal{M}(H;t)$.

Next, if a vertex $v$ is isolated in $G$, i.e., $v$ has no
neighbor, then $E_{v,B(v)}$ is in $\mathcal{L}$ and
$\mathcal{L}=\mathcal{L}^v_0\oplus \langle E_{v,B(v)}\rangle$.
This case is actually studied in the second part, in the proof of
lemma \ref{minor}. In this case, as we already mentioned, for all
$x$, $\mathcal{L}\vert^v_x$ is the same and equal to
$\mathcal{L}\vert^v_0$.

\begin{thm}\label{isoother}{ \begin{equation*}
\mathcal{M}(\mathcal{L}; t)=
\begin{cases} (q-1)t\mathcal{M}(\mathcal{L}\vert^v_0; t) & if\ $v$\ is\ isolated\ in\
$G$,
\\
(q-1)\sum_{x\in \mathbf{K}^*}\mathcal{M}(\mathcal{L}\vert^v_x; t)
& otherwise.
\end{cases}
\end{equation*}

}\end{thm}

\begin{proof}{ Suppose that $v$ is isolated is $G$. Then $\mathcal{L}=\mathcal{L}^v_0\oplus
\langle E_{v,B(v)}\rangle$, and for any complete vector $C\in
\mathbf{K}^V$, we have $$\mathcal{L}\cap
\hat{C}=\big(\mathcal{L}\vert^v_0\cap \hat{C}\vert^v\big)\oplus
\big(\langle E_{v,B(v)}\rangle \cap\langle E_{v, C(v)}\rangle
\big).$$ Therefore,

\begin{eqnarray*} \mathcal{M}(\mathcal{L}; t) & = & \sum_{C\in \mathbf{K}^V}
(t-q)^{\ dim(\mathcal{L}\cap \hat{C})}\\
&=&\sum_{C\in \mathbf{K}^V} (t-q)^{\
dim\big(\mathcal{L}\vert^v_0\cap \hat{C}\vert^v\big)\oplus
\big(\langle E_{v,B(v)}\rangle \cap\langle E_{v, C(v)}\rangle
\big)}\\
&=& \sum_{C\in \mathbf{K}^V} (t-q)^{\ dim\
(\mathcal{L}\vert^v_0\cap \hat{C}\vert^v)}
(t-q)^{\delta_{(B(v)\sim C(v)) }}\\
& =&\sum_{C\in \mathbf{K}^{V-\{v\}}}\big( (q^2-q)+(q-1)(t-q)\big)
(t-q)^{\ dim\ (\mathcal{L}\vert^v_0\cap \hat{C})}\\
&=&(q-1)t\mathcal{M}(\mathcal{L}\vert^v_0; t),
\end{eqnarray*}
where, all summations are over the complete vectors.

Now, assume that $v$ is not isolated and $g_{vw}$ is non-zero for
some $w\in V$. Hence, there exists a vector $C_1\in \mathcal{L}$
such that $C_1(v)\sim A(v)$. Also, we already know that for some
$C_2\in \mathcal{L}$, we have $C_2(v)\sim B(v)$, and $A(v)\nsim
B(v)$. Therefore, the $v$-th coordinates of vectors in
$\mathcal{L}$ cover the whole space $\mathbf{K}$. Thus,

\begin{eqnarray*}\mathcal{M}(\mathcal{L}; t) &=&\sum_{C\in \mathbf{K}^V}
(t-q)^{dim(\mathcal{L}\cap \hat{C})} \\
&=&\sum_{x\in \mathbf{K}^*}\ \,\sum_{C, C(v)\sim x}
(t-q)^{dim(\mathcal{L}\cap \hat{C})}\\
&=&\sum_{x\in \mathbf{K}^*}\ \,\sum_{C, C(v)\sim x}
(t-q)^{dim(\mathcal{L}^v_x\cap \hat{C})}\\
& =&\sum_{x\in \mathbf{K}^*}\ \,\sum_{C, C(v)\sim x}
(t-q)^{dim(\mathcal{L}\vert^v_x\cap \hat{C}\vert^v)}\\
&=&\sum_{x\in \mathbf{K}^*}\ \,\sum_{C\in \mathbf{K}^{V-\{v\}}}
(q-1)(t-q)^{dim(\mathcal{L}\vert^v_x\cap \hat{C})}\\
&= &\sum_{x\in \mathbf{K}^*}
(q-1)\mathcal{M}(\mathcal{L}\vert^v_x; t),
\end{eqnarray*} (once
again, all summations are over the complete vectors).

}\end{proof}

\vspace{3mm}

Consider a graph $G$. As usual, by $G-\{v\}$, we mean the graph
obtained from $G$ by deleting vertex $v$.

\begin{lem}\label{deletion}{ Let $\mathcal{L}$ be an isotropic system with
graphic presentation $(G, A, B)$. Then $G-\{v\}$ is a fundamental
graph of $\mathcal{L}\vert^v_{A(v)}$.

 }\end{lem}

\begin{proof}{ We already know that the rows of $(I\mid G)D(A, B)$ form a basis
for $\mathcal{L}$. The $v$-th coordinate of each row, except the
$v$-th row, is either zero or equivalent to $A(v)$. Therefore, the
rows of $(I\mid G)D(A, B)$, except the $v$-th row form a basis for
$\mathcal{L}^v_{A(v)}$. Thus, by deleting the $v$-th row and the
$v$-th column of $G$, and also the $v$-th coordinates of $A$ and
$B$, we come up with a graphic presentation of
$\mathcal{L}\vert^v_{A(v)}$, meaning that $G-\{v\}$ is a
fundamental graph of $\mathcal{L}\vert^v_{A(v)}$.

}\end{proof}

\begin{thm}\label{operdel}{\hfill

\vspace{3mm}

\noindent{\rm (i)} If $v$ is an isolated vertex in $G$, then
$G-\{v\}$ is a fundamental graph of $\mathcal{L}\vert^v_0$.

\vspace{3mm}

\noindent{\rm (ii)} If $w$ is a neighbor of $v$ in $G$ then $q$
graphs $G*_rv-\{v\}$, $r\in \mathbf{F}_q$, together with
$G*_{-g^{-2}_{vw}}w*_1v-\{v\}$ are fundamental graphs of
$\mathcal{L}\vert^v_x$ for $x\in \mathbf{K}^*$.

}\end{thm}

\begin{proof}{ Part {\rm (i)} is a direct consequece of lemma \ref{deletion}.
To prove {\rm (ii)}, using lemma \ref{deletion}, we should show
that for each $x\in \mathbf{K}^*$, there exists an Eulerian vector
related to one of these graphs, such that its $v$-th coordinate
is $x$.

Suppose that $(G, A, B)$ is a graphic presentation of
$\mathcal{L}$. By theorem \ref{locswit}, for any $r\in
\mathbf{F}_q$, $(G*_rv, A+rB_v, B+rg^2(v)\times A)$ is also a
graphic presentation of $\mathcal{L}$. The $v$-th coordinate of
the Eulerian vector of this presentation is $A(v)+rB(v)$.
Therefore, $G*_rv-\{v\}$ is a fundamental graph of
$\mathcal{L}\vert^v_{A(v)+rB(v)}$.

Now notice that, $\langle A(v), B(v)\rangle\neq 0$. Therefore, for
$r$ varies in $\mathbf{F}_q$, $A(v)+rB(v)$ are all different
elements of $\mathbf{K}^*$, and there are $q$ of them. The only
element of $\mathbf{K}^*$ not obtained in this way is $B(v)$.
Consider the fundamental graph $G*_sw*_1v$, where
$s=-g^{-2}_{vw}$. By theorem \ref{locswit},
$A+sB_w+(B_v+sg^2_{vw}A_v)$ is an Eulerian vector for this
fundamental graph, and the $v$-th coordinate of this vector is
$(sg^2_{vw}+1)A(v)+B(v)=B(v)$. Thus, $G*_sw*_1v-\{v\}$ is a
fundamental graph of $\mathcal{L}\vert^v_{B(v)}$.

 }\end{proof}

\begin{cor}\label{recmar}{
If $v$ is isolated in $G$ then $$\mathcal{M}(G; t)=
(q-1)t\mathcal{M}(G-\{v\}; t),$$ otherwise, if $w$ is a neighbor
of $v$ then $$\mathcal{M}(G;
t)=(q-1)\Big[\mathcal{M}(G*_{-g^{-2}_{vw}}w*_1v; t)+\sum_{r\in
\mathbf{F}_q}\mathcal{M}(G*_rv-\{v\}; t)\Big].$$

}\end{cor}

\subsection{Estimation of $\epsilon(G)$}

The final step to evaluate $\epsilon(G)$ is the following one,
which is valid due to the fact that all of the graphs described
in the right hand side of the formula in corollary \ref{recmar}
have $n-1$ vertices. Indeed, one can observe that the number of
graphs in the right hand side is equal to $q+1$, and hence,
$$\max_{G:\ |G|=n} |\mathcal{M}(G;t)| \leq (q-1)\cdot \max\{t,q+1\}
\cdot \max_{H:\ |H|=n-1} |\mathcal{M}(H;t)|.$$ Moreover, when the
graph $G$ has only one vertex, i.e., $n=1$, we have
$$|\mathcal{M}(G;t)| \leq (q^2-1)\cdot \max \{|t-q|,1\}.$$
Putting together these two statements, we conclude the following
corollary.

\begin{cor}\label{recmar2}{ For a graph $G$ with $n$ vertices, the
Tutte-Martin polynomial $\mathcal{M}(G;t)$ can be estimated as
follows:
$$|\mathcal{M}(G;t)| \leq (q^2-1) \cdot [(q-1)\cdot \max\{t,q+1\}]^{(n-1)} \cdot \max
\{|t-q|,1\},$$ and by setting $t=q$, we obtain that:
$$\epsilon(G) \leq (q^2-1)^n.$$

}\end{cor}

\vspace{3mm}

We can even derive a lower bound for $\epsilon(G)$ as well.
$$\min_{G:|G|=n} \mathcal{M}(G;q) \geq (q-1)\cdot \min\{t,q+1\} \cdot \min_{H:\
|H|=n-1} |\mathcal{M}(H;q)|.$$ On the other hand, when the graph
has just one vertex, $\epsilon(G) \geq 1$, and hence in general,

$$\epsilon(G) \geq (q^2-q)^{n-1}.$$
Thus, the proof of the following theorem is now complete.

\begin{thm}\label{bound}{ The number of Eulerian vectors for a graph
(or equivalently for an isotropic system) with $n$ vertices
satisfies the following property:

$$(q^2-q)^{n-1} \leq \epsilon(G) \leq (q^2-1)^n.$$
In particular, when the graph is ordinary (binary), i.e., when
$q=2$, we have:
$$2^{n-1} \leq \epsilon(G) \leq 3^n.$$
}\end{thm}
 \hfill{$\Box$}

\section{The number of classes of local equivalency}

As mentioned earlier, we use the formula given in corollary
\ref{kerel2}, as well as the estimations for the number of
Eulerian vectors given in the previous section, in order to give
a bound for $l(G)$, the number of graphs locally equivalent to
$G$.

By corollary \ref{kerel2}, if $G$ is connected, then
$$l(G) \leq \frac{(q-1)\epsilon(\mathcal{L})}{\lambda(\mathcal{L})} \leq (q-1)\epsilon(\mathcal{L})$$
Taking into account the estimation of
$\epsilon(\mathcal{L})=\epsilon(G)$ presented in the theorem
\ref{bound}, we come up with an upper bound for $l(G)$, given
that the number of graphs with $n$ vertices is exactly
$q^{\frac{n^2}{2}-\frac{n}{2}}$.

\begin{thm}\label{main2}{ \hfill

\vspace{3mm}

\noindent{\rm (i)} The number of graphs locally equivalent to a
connected graph is at most $(q-1)(q^2-1)^n$ which is bounded
above by $q^{2n+1}$, $n$ being the number of vertices of the
graph.

\vspace{3mm}

\noindent{\rm (ii)} $\mathcal{C}(n)$, the number of classes of
local equivalency of connected graphs with $n$ vertices satisfies:
$$q^{\frac{n^2}{2}-\frac{5n}{2}-1} \leq \mathcal{C}(n) \leq q^{\frac{n^2}{2}-\frac{n}{2}}. $$
In other words, $$\mathcal{C}(n) = q^{\frac{n^2}{2}-O(n)}.$$ In
particular, for the usual (binary) graphs, i.e., when $q=2$, the
number of graphs locally equivalent to a graph is at most $3^n$
and the number of classes of local equivalency is
$$\mathcal{C}(n) = 2^{\frac{n^2}{2}-O(n)}.$$

\hfill$\Box$

}\end{thm}

\section{Conclusion}

We developed a method to compute number of graphs locally
equivalent to a given one. Using this method, we bounded this
number for an arbitrary graph. Also, we found an approximation of
the number of equivalency classes. That is, $\mathcal{C}(n) =
q^{\frac{n^2}{2}-O(n)}.$ Notice that, number of all graphs is
$q^{\frac{n^2}{2}-\frac{n}{2}}$. Therefore, this estimation says
that number of equivalency classes is almost the same as number
of all graphs.

We got these results by developing a reach theory of isotropic
systems and also local operation over graphs, and it seems that
this theory can be used for other problems in this area, too.

\vspace{3mm}

\noindent {\bf Acknowledgement.} The authors a greatly thankful
to Peter Shor for all his support and helpful advice.

\end{document}